\documentclass[11pt, reqno]{amsart}
\usepackage{mathrsfs}
\usepackage[active]{srcltx}
\usepackage{mathrsfs}
\usepackage{color}

\usepackage{amsthm}
\usepackage{amssymb}

\usepackage{color}
\usepackage[unicode]{hyperref}
\hypersetup{
	colorlinks = true,%
	citecolor = [rgb]{0.0,0.0,0.9},
	filecolor=black,%
	linkcolor = [rgb]{0.65,0.0,0.0},%
	anchorcolor = red,
	pagecolor = red,
	urlcolor= [rgb]{0.65,0.0,0.0}
}

\newtheorem{proposition}{Proposition}[section]
\newtheorem{theorem}{Theorem}[section]

\newtheorem{remark}{Remark}[section]

\numberwithin{equation}{section}

\begin{document}
	
	\title[Interpolation between Brezis-V\'azquez and Poincar\'e inequalities]
	{Interpolation between Brezis-V\'azquez and Poincar\'e inequalities on nonnegatively curved spaces:\\ sharpness and rigidities}
	
\vspace{-0.9cm}
	\author{Alexandru Krist\'aly}
\address{Institute of Applied Mathematics, \'Obuda
	University, 
	Budapest, Hungary \& Department of Economics, Babe\c
	s-Bolyai University, Cluj-Napoca, Romania}
	\email{kristaly.alexandru@nik.uni-obuda.hu; alexandrukristaly@yahoo.com 
}
	
	\author{Anik\'o Szak\'al}
		\address{Institute of Applied Mathematics, \'Obuda
		University, 
		Budapest, Hungary}
	\email{szakal@uni-obuda.hu}

	


\newcommand{\abs}[1]{\lvert#1\rvert}

\newcommand{\blankbox}[2]{%
	\parbox{\columnwidth}{\centering
		\setlength{\fboxsep}{0pt}%
		\fbox{\raisebox{0pt}[#2]{\hspace{#1}}}%
	}%
}

\vspace{-0.6cm}

\begin{abstract}
{\footnotesize \noindent This paper is devoted to investigate an  interpolation inequality between the Brezis-V\'azquez and Poincar\'e inequalities (shortly, BPV inequality) on  nonnegatively curved spaces. As a model case, we first prove that the BPV inequality holds on any Minkowski space, by fully characterizing the existence and shape of its extremals. We then prove that if a complete Finsler manifold with nonnegative Ricci curvature supports the BPV inequality, then its flag curvature is identically zero. In particular, we deduce that  a Berwald space of nonnegative Ricci curvature supports the BPV inequality if and only if it is isometric to a Minkowski space.  Our arguments explore fine properties of Bessel functions, comparison principles, and anisotropic symmetrization on Minkowski spaces. As an application, we characterize the existence of nonzero solutions for a quasilinear PDE involving the Finsler-Laplace operator and a Hardy-type singularity on Minkowski spaces where the sharp BPV inequality plays a crucial role. The results are also new  in the Riemannian/Euclidean setting. 	 
	     }
\end{abstract}

\maketitle

\vspace{-0.9cm}

\section{Introduction}

One of the most spectacular improvements of the classical unipolar Hardy inequality  is due to Brezis and V\'azquez \cite{BV} by establishing that for every bounded domain $\Omega\subset \mathbb R^n$ $(n\geq 2)$ with $0\in \Omega$ one has 
$$\int_\Omega |\nabla u(x)|^2{\rm d}x\geq \frac{(n-2)^2}{4}\int_\Omega \frac{u(x)^2}{|x|^2}{\rm d}x+j_0^2\left(\frac{\omega_n}{|\Omega|}\right)^\frac{2}{n} \int_\Omega u(x)^2{\rm d}x,\ \ \forall u\in W_0^{1,2}(\Omega),\eqno{(\textbf{BV})} $$
where $j_0=2.4048$ is the first positive zero of the Bessel function of first kind $J_0$, while $|\Omega|$ and $\omega_n$ denote the volumes of the set $\Omega$ and the  $n$-dimensional Euclidean unit ball, respectively. In the limit case $n=2$, the inequality $(\textbf{BV})$ reduces precisely to the optimal Poincar\'e inequality.

The aforementioned inequalities constitute a continuous source of inspiration for further investigations not only in the Euclidean setting, see e.g. Adimurthi,  Chaudhuri and Ramaswamy \cite{ACR}, Barbatis,  Filippas and Tertikas \cite{Barbatis},
Ghoussoub and Moradifam \cite{GM, GM-2}, but also on curved spaces. More precisely, such Sobolev-type inequalities behave quite naturally on \textit{Hadamard manifolds} (simply connected, complete Riemannian/Finsler manifolds with nonpositive sectional/flag curvature), as shown e.g. by  Carron \cite{Caron, Carron-2}, Berchio,  Ganguly and Grillo \cite{Grillo}, D'Ambrosio and  Dipierro \cite{D-D-olaszok}, Kombe and \"Ozaydin \cite{KO-1, KO-2}, 
Farkas, Krist\'aly and Varga \cite{FKV}, Krist\'aly \cite{Kristaly-JMPA}, Yang, Su and  Kong \cite{YSK}. This fact is not surprising since Hadamard manifolds are diffeomorphic to Euclidean spaces.    

Our paper is devoted to study an  inequality  on \textit{nonnegatively curved spaces} whose limit cases are the Brezis-V\'azquez  and (not necessarily the 2-dimensional) Poincar\'e inequalities.  

In order to formulate the interpolation inequality,  
let $({M},F)$ be a complete $n$-dimensional reversible  Finsler manifold  ($n\geq 2$ be an integer) and $\alpha\in \left[0,\frac{n-2}{2}\right]$ be fixed.  If $\Omega\subset M$ is a bounded open set and $x_0\in \Omega$, 
we consider the \textit{Brezis-Poincar\'e-V\'azquez inequality} 
$$\int_{\Omega} F_*(x,Du(x))^2 {\rm d}V_F(x) \geq
\left[\frac{(n-2)^2}{4}-\alpha^2\right]\int_{\Omega}\frac{u(x)^2}{{d}_F(x_0,x)^2}{\rm d}V_F(x)+S_{\alpha}(\Omega)\int_{\Omega}
{u(x)^2}{\rm d}V_F(x),\ \ \forall u\in C_0^\infty(\Omega),
\eqno{(\textbf{BPV})}$$
where  
$$S_{\alpha}(\Omega):=j_{\alpha}^2\left(\frac{\omega_n}{{\rm Vol}_F(\Omega)}\right)^\frac{2}{n},$$ and $j_{\alpha}$ is the first
positive zero of the Bessel function of the first kind $J_\alpha$. 
Hereafter, $F_*$, $d_F$,   ${\rm d}V_F$ and Vol$_F$ denote the polar transform, the metric function, the canonical measure and Finslerian volume on $(M,F)$, respectively; for details, see Section \ref{section-pre}.


In the classical Euclidean setting,    $(\textbf{BPV})$ reduces to the \textit{Brezis-V\'azquez inequality} $(\textbf{BV})$ when $\alpha=0$, and to the optimal \textit{Poincar\'e inequality} when  $\alpha= \frac{n-2}{2}$. 

In fact, our first main result shows that 
$(\textbf{BPV})$ holds on  \textit{Minkowski spaces}, the simplest Finslerian structures with vanishing flag curvature (i.e., $\mathbb R^n$ endowed with an arbitrary smooth norm). Without  loss of generality, the Minkowski norm in $(\mathbb R^n,F)$ is scaled such that  $B_0^F(1)=\{x\in \mathbb R^n:F(x)<1\}$ has volume $\omega_n.$ As usual, a set $\Omega\subset \mathbb R^n$ has a Wulff shape if it is homothetic to  $B^F_0(1)$.  For further use, let $l_F\in (0,1]$ be the uniformity constant associated with $F$; we note that  $l_F=1$ if and only if $F$ is Euclidean, see Section \ref{section-pre}. 
By using anisotropic symmetrization arguments, see Alvino,  Ferone,  Lions and Trombetti \cite{AIHP-Lions} and Van Schaftingen \cite{vanSch}, and fine convexity properties of the Hardy functional involving the uniformity constant $l_F$ on $(\mathbb R^n,F)$, we prove the following result.  

\begin{theorem}\label{theorem-extremals}
	Let $(\mathbb R^n,F)$ be a Minkowski space, $n\geq 2,$ and fix  $\alpha\in \left[\frac{n-2}{2}\sqrt{1-l_F^2},\frac{n-2}{2}\right]$. Then inequality {\rm $(\textbf{BPV})$} holds for every open bounded set $\Omega\subset \mathbb R^n$ and $x_0\in \Omega$. 
	
	Moreover, given an open set $\Omega\subset \mathbb R^n$,   equality holds in {\rm $(\textbf{BPV})$}  for some function belonging to the Sobolev space $W_0^{1,2}(\Omega)$ if and only if $\Omega$ has a Wulff shape and either $\alpha=0$ when $n=2$, or $\alpha>0$ when $n\geq 3;$  
	in such cases, the extremal function has the form 
	\begin{equation}\label{u-csillag-megoldas}
	u^\star(x)=F(x)^\frac{2-n}{2}J_\alpha\left({\sqrt{S_{\alpha}(\Omega)}}F(x)\right),\ x\in \Omega^\star, 
	\end{equation}
	where $\Omega^\star$ is the anisotropic symmetrization of $\Omega.$
\end{theorem}

Having the flat case (Theorem \ref{theorem-extremals}), a natural question arises: what about the  $(\bf{BPV})$ inequality on nonnegatively curved  Finsler manifolds? The answer is given in the following rigidity result. 


\begin{theorem}\label{rigidity}  Let $(M,F)$ be a complete $n$-dimensional reversible Finsler manifold $(n\geq 2)$ with nonnegative $n$-Ricci curvature, 
	and  $\alpha\in \left[0,\frac{n-2}{2}\right]$ be such that 
	$\alpha>0$ whenever $n\geq 3.$ If {\rm $(\textbf{BPV})$} holds for every 
	open bounded
	set $\Omega\subset M$ and $x_0\in \Omega$, then  the flag curvature of $(M,F)$ is identically zero.  
\end{theorem}

The proof of Theorem \ref{rigidity} requires a fine analysis of Bessel functions combined with the Bishop-Gromov volume comparison principle  on Finsler manifolds. 

 Theorem \ref{rigidity} is new in the Riemannian setting as well; however, its conclusion in the  particular case $\alpha=\frac{n-2}{2}$ can be obtained by the Rayleigh-Faber-Krahn inequality established by  Cheng \cite{Cheng}. Indeed, when $(M,g)$ is an $n$-dimensional Riemannian manifold with nonnegative Ricci curvature endowed with its natural metric $d_g$ and canonical measure ${\rm d}V_g$, we have the \textit{Rayleigh-Faber-Krahn inequality} 
 \begin{equation}\label{Cheng-res}
 \mu_1(B_{x}(\rho)):=\inf_{u\in W_0^{1,2}(B_{x}(\rho))\setminus \{0\}}\frac{\displaystyle\int_{B_{x}(\rho)} |Du|(x)^2{\rm d}V_g(x)}{\displaystyle\int_{B_{x}(\rho)} u(x)^2{\rm d}V_g(x)}\leq \mu_1(B^e_{0}(\rho))=\frac{j_{n/2-1}^2}{\rho^2},
 \end{equation}
 where $B_{x}(\rho)=\{y\in M:d_g(x,y)<\rho\}$ and $B^e_{0}(\rho)$ is the $n$-dimensional Euclidean ball with center $0$ and radius $\rho>0;$ moreover,  equality holds in (\ref{Cheng-res}) if and only if $B_{x}(\rho)$ is isometric to  $B^e_{0}(\rho)$,  see Cheng \cite{Cheng}. In this Riemannian setting, the validity of the inequality $(\bf{BPV})$ with  $\alpha=\frac{n-2}{2}$ (i.e., Poincar\'e inequality) implies equality in (\ref{Cheng-res}), thus the conclusion in Theorem \ref{rigidity} directly follows by Cheng's result. However, Cheng's approach -- based on a careful analysis of Jacobi fields on normal coordinates of $(M,g)$ -- cannot be adapted to our  setting, where some singular terms also occur when $\alpha\neq \frac{n-2}{2}$.

Theorems \ref{theorem-extremals} and \ref{rigidity} can be elegantly  summarized on Berwald spaces,  by providing an analytic characterization of Minkowski spaces through the $(\textbf{BPV})$ inequality. 

\begin{theorem}\label{rigidity-Minkowski}  Let $(M,F)$ be a complete $n$-dimensional reversible Berwald space $(n\geq 2)$ having nonnegative Ricci curvature,  uniformity constant $l_F\in (0,1]$, 
	and fix  $\alpha\in \left[\frac{n-2}{2}\sqrt{1-l_F^2},\frac{n-2}{2}\right]$ 
	such that 
	$\alpha>0$ whenever $n\geq 3.$
	Then the following two statements are equivalent$:$
	\begin{itemize}
		\item[(i)] {\rm $(\textbf{BPV})$} holds for every 
		open bounded
		set $\Omega\subset M$ and $x_0\in \Omega;$
		\item[(ii)] $(M,F)$ is isometric to a Minkowski space. 
	\end{itemize}
\end{theorem}

As an application of the $(\textbf{BPV})$ inequality, we consider on a Minkowski space $(\mathbb R^n,F)$ the following quasilinear Dirichlet problem 
\[ \   \left\{ \begin{array}{lll}
-\Delta_F u(x)-\left[\frac{(n-2)^2}{4}-\alpha^2\right]\frac{u(x)}{F(x)^2}+\lambda u(x)  = |u(x)|^{p-2 }u(x),& &  x\in B_0^F(1); \\
u\geq 0,\ u\in W_0^{1,2}(B_0^F(1)),
\end{array}\right. \eqno{({\mathcal P}_{\alpha,\lambda})}\]
where $\Delta_F u={\rm div}(J^*(Du(x)))$ is the Finsler-Laplace operator on $(\mathbb R^n,F)$, $J^*$ being the Legendre transform associated to $F$, see Section \ref{section-pre}. The following result characterizes the existence of nonzero solutions of problem $({\mathcal P}_{\alpha,\lambda})$ depending on the parameters $\alpha,\lambda\in \mathbb R.$ As usual, $2^*$ denotes the critical Sobolev exponent ($2^*=2n/(n-2)$ if $n\geq 3$ and $2^*=\infty$ if $n=2$).

\begin{theorem}\label{alkalmazas}
Let $(\mathbb R^n,F)$ be a Minkowski space, $n\geq 2,$ and fix  $\alpha\in \left[\frac{n-2}{2}\sqrt{1-l_F^2},\frac{n-2}{2}\right]$ such that  $\alpha>0$ whenever $n\geq 3.$ Let $p\in (2,2^*)$ be fixed. Then problem $({\mathcal P}_{\alpha,\lambda})$ has a nonzero solution if and only if $\lambda>-j_\alpha^2.$
\end{theorem}

Theorem \ref{alkalmazas} is known in the special case when $F$ is Euclidean and $\alpha=\frac{n-2}{2}$ (thus, the singular term disappears), see Willem \cite{Willem}. The proof of Theorem \ref{alkalmazas} is variational, based on the  mountain pass theorem and the validity of the sharp inequality $(\textbf{BPV})$ on Minkowski spaces. \\

%
%
%
%
%

The organization of the paper is the following. 
In Section \ref{section-pre} we recall basic notions from Finsler geometry (flag curvature, Ricci curvature, Bishop-Gromov volume comparison principle). In Section \ref{section-3}, before presenting the proof of Theorem \ref{rigidity}, we recall some basic results from the theory of anisotropic symmetrization on Minkowski spaces.  In Section \ref{section-4} we prove Theorems \ref{rigidity} and \ref{rigidity-Minkowski}; to complete this,  we first establish some properties of Bessel functions which are interesting in their own right. Finally, in Section \ref{section-5} we prove Theorem \ref{alkalmazas}.

\section{Preliminaries on Finsler manifolds}\label{section-pre}

Let $M$ be a connected $n$-dimensional $C^{\infty}$-manifold and $TM=\bigcup_{x \in M}T_{x} M $ be its tangent bundle. 
The pair $(M,F)$ is called a reversible {\it Finsler manifold} if the
continuous function $F:TM\to [0,\infty)$ satisfies the
conditions:
\begin{itemize}
	\item[{\rm (a)}] $F\in C^{\infty}(TM\setminus\{ 0 \});$
	
	\item[{\rm (b)}] $F(x,tv)=|t|F(x,v)$ for all $t\in \mathbb R$ and $(x,v)\in TM;$
	
	\item[{\rm (c)}] the $n \times n$ matrix
	\begin{equation}\label{g_v}
	g_{(x,v)}:=[g_{ij}(x,v)]_{i,j=1,...,n}=\left[\frac{1}{2}\frac{\partial^2 }{\partial v^{i}
		\partial v^{j}}F^2(x,v)\right]_{i,j=1,...,n},
	\quad \text{where}\ v=\sum_{i=1}^n v^i \frac{\partial}{\partial x^i},
	\end{equation}
	is positive definite for all $(x,v)\in TM\setminus\{ 0 \}$. We will
	denote by $g_v$ the inner product on $T_xM$ induced by
	\eqref{g_v}. 
\end{itemize}
If $g_{ij}(x)=g_{ij}(x,v)$ is
independent of $v$ then $(M,F)=(M,g)$ is called a \textit{Riemannian
	manifold}. A \textit{Minkowski space} consists
of a finite dimensional vector space $V$ (identified with $\mathbb{R%
}^{n}$) and a Minkowski norm which induces a Finsler metric on $V$
by translation, i.e., $F(x,v)$ is independent on the base point $x$;
in such cases we often write $F(v)$ instead of $F(x,v)$. A Finsler
manifold $(M,F)$ is called a {\it locally Minkowski space} if any
point in $M$ admits a local coordinate system $(x^i)$ on its
neighborhood such that $F(x,v)$ depends only on $v$ and not on $x$.

For every
$(x,\alpha)\in T^*M$, the \textit{polar transform} (or, co-metric)
of $F$ is given  by
\begin{equation}  \label{polar-transform}
F_*(x,\alpha)=\sup_{v\in T_xM\setminus
	\{0\}}\frac{\alpha(v)}{F(x,v)}.
\end{equation}
Note that for every $x\in M$, the function $F_*(x,\cdot)$ is a
Minkowski
norm on $T_x^*M.$

The number
\begin{equation*}
l_{F}=\inf_{x\in M}l_F(x),\ \ \ \mathrm{where}\ \ \ l_F(x):=
\inf_{y,v,w\in T_xM\setminus
	\{0\}}\frac{g_{(x,v)}(y,y)}{g_{(x,w)}(y,y)},
\end{equation*}
is the \textit{uniformity constant} associated with $F$ which measures how far $(M,F)$ and $%
(M,F_*)$ are from Riemannian structures. In fact, one can see that
$l_{F}\leq 1$, and $l_{F}= 1$ if and only if $(M,F)$ is a Riemannian
manifold. When $(\mathbb R^n,F)$ is a Minkowski space, we have that $l_F\in (0,1]$. The definition of $l_{F}$ in turn shows that
\begin{equation}  \label{eq:2uni}
\left[F_{*}( x,{t\alpha+(1-t)\beta} )\right]^2 \le t\left[F_{*}(x,\alpha)\right]^2
+(1-t)\left[F_{*}(x,\beta)\right]^2 -{l_{F}}t(1-t) \left[F_{*}(x,\beta-\alpha)\right]^2
\end{equation}
for all $x\in M$, $\alpha,\beta \in T_x^*M$ and $t\in [0,1]$. 

Let $\pi ^{*}TM$ be the pull-back bundle of the tangent bundle $TM$
generated by the natural projection $\pi:TM\setminus\{ 0 \}\to M,$
see Bao, Chern and Shen \cite{BCS}. The vectors of the pull-back
bundle $\pi
^{*}TM$ are denoted by $(v;w)$ with $(x,y)=v\in TM\setminus\{ 0 \}$ and $%
w\in T_xM.$ For simplicity, let $\partial_i|_v=(v;\partial/\partial
x^i|_x)$ be the natural local basis for $\pi ^{*}TM$, where $v\in
T_xM.$ One can introduce on $\pi ^{*}TM$ the \textit{fundamental
	tensor} $g$
by
$
g_{(x,v)}:=g_v=g(\partial_i|_v,\partial_j|_v)=g_{ij}(x,y),
$
where $v=y^i{(\partial}/{\partial x^i})|_x,$ see (\ref{g_v}). Unlike the Levi-Civita
connection \index{Levi-Civita connection} in the Riemannian case,
there is no unique natural connection in the Finsler geometry. Among
these connections on the pull-back bundle $\pi ^{*}TM,$ we choose a
torsion-free and almost
metric-compatible linear connection on $\pi ^{*}TM$, the so-called \textit{%
	Chern connection}. The coefficients of the Chern connection are
denoted by $\Gamma_{jk}^{i}$, which are instead of the well-known
Christoffel symbols from Riemannian geometry.

A Finsler manifold is of \textit{Berwald type} if the coefficients $%
\Gamma_{ij}^{k}(x,y)$ in natural coordinates are independent of $y$.
It is clear that Riemannian manifolds and $($locally$)$ Minkowski
spaces are
Berwald spaces. The Chern connection induces on $\pi ^{*}TM$ the \textit{%
	curvature tensor} $R$.  The Finsler manifold is
{\it complete} if every
geodesic segment $\sigma:[0,a]\to M$ can be extended to $\mathbb R$.

Let $u,v\in T_xM$ be two non-collinear vectors and $\mathcal{S}=\mathrm{span}%
\{u,v\}\subset T_xM$. By means of the curvature tensor $R$, the
\textit{flag curvature} associated with the flag $\{\mathcal{S},v\}$ is 
\begin{equation}  \label{ref-flag}
\mathbf{K}(\mathcal{S};v) =%
\frac{g_v(R(U,V)V, U)}{g_v(V,V) g_v(U,U) - g_v^{2}(U,V)},
\end{equation}
where $U=(v;u),V=(v;v)\in \pi^*TM.$ If $(M,F)$ is Riemannian, the
flag
curvature reduces to the sectional curvature which depends only on $\mathcal{S}$.  

Take $v\in T_xM$ with $F(x,v)=1$ and
let $\{e_i\}_{i=1}^n$ with $e_n=v$ be an orthonormal basis of
$(T_xM,g_v)$ for $g_v$ from \eqref{g_v}. Let $\mathcal S_i={\rm
	span}\{e_i,v\}$ for $i=1,...,n-1$. Then the {\it Ricci curvature} of
$v$ is defined by ${\rm Ric}(v):=\sum_{i=1}^{n-1}\textbf{K}(\mathcal S_i;v)$.

Let $\mu$ be a positive smooth measure on $(M,F)$. Given $v \in
T_xM \setminus \{0\}$, let $\sigma:(-\varepsilon,\varepsilon)
\to M$ be the geodesic with $\dot{\sigma}(0)=v$ and
decompose $\mu$ along $\sigma$ as
$\mu=e^{-\psi}\mathrm{vol}_{\dot{\sigma}}$, where
$\mathrm{vol}_{\dot{\sigma}}$ denotes the volume form of the
Riemannian structure $g_{\dot{\sigma}}$. For $N \in
[n,\infty]$, the \emph{$N$-Ricci curvature} $\mathrm{Ric}_N$ is
defined by
\[ \mathrm{Ric}_N(v):=\mathrm{Ric}(v) +(\psi \circ \sigma)''(0)
-\frac{(\psi \circ \sigma)'(0)^2}{N-n}, \]
where the third term is understood as $0$ if $N=\infty$ or if $N=n$
with $(\psi \circ \sigma)'(0)=0$, and as $-\infty$ if $N=n$ with
$(\psi \circ \sigma)'(0) \neq 0$.

Let $\sigma: [0,r]\to M$ be a piecewise smooth curve. The value $%
L_F(\sigma)= \displaystyle\int_{0}^{r} F(\sigma(t), \dot\sigma(t))\,{\text d}%
t $ denotes the \textit{integral length} of $\sigma.$ For
$x_1,x_2\in M$,
denote by $\Lambda(x_1,x_2)$ the set of all piecewise $C^{\infty}$ curves $%
\sigma:[0,r]\to M$ such that $\sigma(0)=x_1$ and $\sigma(r)=x_2$.
Define the \textit{metric function} $d_{F}: M\times M
\to[0,\infty)$ by
\begin{equation}  \label{quasi-metric}
d_{F}(x_1,x_2) = \inf_{\sigma\in\Lambda(x_1,x_2)} L_F(\sigma).
\end{equation}
The metric ball with center $x\in M$ and radius $\rho>0$ is defined by $B_x(\rho)=\{y\in M:d_F(x,y)<\rho\}.$

Let $\{{\partial}/{\partial x^i} \}_{i=1,...,n}$ be a local basis
for the tangent bundle $TM,$ and $\{\mathrm{d}x^i \}_{i=1,...,n}$ be
its dual basis
for $T^*M.$ Consider $\tilde B_x(1)=\{y=(y^i):F(x,y^i \partial/\partial x^i)< 1\}\subset \mathbb R^n$. The \textit{Busemann-Hausdorff volume form}  is defined by
\begin{equation}  \label{volume-form}
{\text d}V_F(x)=\sigma_F(x){\text
	d}x^1\wedge...\wedge {\text d}x^n,
\end{equation}
where $\sigma_F(x)=\frac{\omega_n}{|\tilde B_x(1)|}$. The \textit{%
	Finslerian volume} of an open set $S\subset M$ is
Vol$_F(S)=\displaystyle\int_S {\text d}{V}_F(x)$. When $(\mathbb{R}^n,F)$ is a Minkowski space, then 
$d_F(x_1,x_2)=F(x_2-x_1)$ and 
on account of (\ref%
{volume-form}), Vol$_F(B_x(\rho))=\omega_n\rho ^n$ for every $\rho>0$ and $%
x\in \mathbb{R}^n$.

%

On any Finsler manifold $(M,F) $ we have  for every $x\in M$ that
\begin{equation}  \label{volume-comp-nullaban}
\lim_{\rho\to 0^+}\frac{\mathrm{Vol}_F(B_x(\rho))}{\omega_n\rho ^n}%
=1.
\end{equation}

Let $(M,F)$ be a complete $n$-dimensional Finsler manifold with
nonnegative $N$-Ricci curvature. Then the Bishop-Gromov volume comparison principle provides that the function  $$\rho\mapsto \frac{{\rm Vol}_F(B_x(\rho))}{\rho^N},\ \ \rho>0,$$ is non-increasing for every $x\in M$. In particular, if $N=n$, then	\begin{equation}\label{volume-comp-altalanos-2}
{{\rm Vol}_F(B_x(\rho))}\leq \omega_n \rho^n,\ \forall x\in M, \ \rho>0.
\end{equation}
Moreover, if equality holds in {\rm (\ref{volume-comp-altalanos-2})}, then the
flag curvature is identically zero, see Ohta \cite{Ohta-Calculus}, Shen \cite{Shen-volume}.

 The (distributional) \textit{derivative} of $%
u:M\to \mathbb R$ at $x\in M$ is 
$Du(x)=\sum_{i=1}^n \frac{\partial }{\partial x^i}u(x)\mathrm{d}x^i, $
and 
due to Ohta and Sturm \cite{Ohta-Sturm}, one has the eikonal equation  
{\
	\begin{equation}  \label{tavolsag-derivalt}
	F_*\left(x,D d_F(x_0,\cdot)(x)\right)=1\ \mathrm{for\ a.e.}\ x\in M.
	\end{equation}%
} 
The {\it Legendre transform}
$J^*:T^*M\to TM$ associates to each element $\xi\in T_x^*M$ the
unique maximizer on $T_xM$ of the map $y\mapsto
\xi(y)-\frac{1}{2}F(x,y)^2$. The {\it gradient} of $u$ is defined by
${\nabla}_F u(x)=J^*(x,Du(x)).$
The {\it Finsler-Laplace operator} is given by 
$${\Delta}_F u={\rm div}({\nabla}_F u),$$ where 
div$(X)=\frac{1}{\sigma_F}\frac{\partial}{\partial
	x^i}(\sigma_F X^i)$ for some vector field $X$  on $M$, and $\sigma_F$ comes from (\ref{volume-form}). 

Consider the \textit{Sobolev space} 
\begin{equation*}
W^{1,2}(M,F):=\left\{ u\in W_{\mathrm{loc}}^{1,2}(M):\displaystyle%
\int_{M}F_{\ast }(x,Du(x))^2{\text d}V_F(x)<+\infty \right\}
,
\end{equation*}%
associated with $(M,F)$
and let $W_{0}^{1,2}(M,F,\mathsf{m})$ be the closure of $C_{0}^{\infty
}(M)$ with respect to the norm
$$\Vert u\Vert _{F}:=\left( \displaystyle\int_{M}F_{\ast }(x,Du(x))^2{\text d}V_F(x)+\displaystyle\int_{M}u(x)^{2}{\text d}V_F(x)\right)
^{1/2}.$$

When $(\mathbb{R}^n,F)$ is a Minkowski space, then $W_{0}^{1,2}(\Omega,F)$ is the usual Sobolev space $W_{0}^{1,2}(\Omega)$ for every open set $\Omega\subset \mathbb R^n$, see Krist\'aly and Ohta \cite{Kri-Ohta}; indeed, in this case there exits $C_0\geq 1$ such that $C_0^{-1}|x|\leq F(x)\leq C_0|x|$ for every $x\in \mathbb R^n.$

\section{Proof of Theorem \ref{theorem-extremals}}\label{section-3}

Before to present the proof of Theorem \ref{theorem-extremals}, we recall some notions and results established in Alvino,  Ferone,  Lions and Trombetti \cite{AIHP-Lions} and Van Schaftingen \cite{vanSch} concerning anisotropic symmetrization. 

Let $(\mathbb R^n,F)$ be a Minkowski space, $n\geq 2$. If $\Omega\subset \mathbb
R^n$ is a measurable set, we denote by $\Omega^\star$ its anisotropic symmetrization defined as the open ball with center $0$ such that $|\Omega|=|\Omega^\star|$. It is clear that $\Omega^*$ has a Wulff shape, homothetic to $B^F_0(1)=\{x\in \mathbb R^n:F(x)<1\}$. 
If $u:\mathbb R^n\to [0,\infty)$ is  a
function, then  $$u^\star(x)=\sup\{c\in \mathbb R:x\in
\{u>c\}^\star\}$$ is the {\it anisotropic $($decreasing$)$
	symmetrization} of $u$. Here, $\{u>c\}=\{x\in \mathbb R^n:u(x)>c\}.$
The following results are valid: 

\begin{itemize}
	\item[$\bullet$] {\it Anisotropic Cavalieri principle} (see 
	\cite[Proposition 2.28]{vanSch}). Let $u:\mathbb
	R^n\to [0,\infty)$ be a function vanishing at infinity with respect
	to $\cdot^\star.$  Then
	$$\int_{\mathbb R^n} u(x)^2{\rm d}x=\int_{\mathbb R^n} u^\star(x)^2{\rm d}x.$$
	\item[$\bullet$] {\it Anisotropic P\'olya-Szeg\H{o} inequality} (see \cite[Theorem 3.1]{AIHP-Lions} and \cite[Theorem 6.8]{vanSch}). If $\Omega\subset \mathbb R^n$ is an open
	set and $u\in W_0^{1,2}(\Omega)_+=\{u\in W_0^{1,2}(\Omega):u\geq 0 \},$ then $u^\star\in
	W_0^{1,2}(\Omega^\star)_+$ and
	$$\int_{\Omega^\star} F_*(D u^\star(x))^2{\rm d}x\leq
	\int_{\Omega} F_*(D u(x))^2{\rm d}x.$$
	\item[$\bullet$] {\it Anisotropic Hardy-Littlewood inequality} (see \cite[Proposition 2.28]{vanSch}). For every open $\Omega\subset \mathbb R^n$  and $u\in W_0^{1,2}(\Omega)_+$, one has 
	$$\int_{\Omega}\frac{u(x)^2}{F(x)^2}{\rm d}x\leq \int_{\Omega^\star}\frac{u^\star(x)^2}{F(x)^2}{\rm d}x.$$
\end{itemize}

\noindent Finally, we recall the \textit{anisotropic Hardy inequality} from \cite{vanSch}: when $n\geq 3$,  for every open set $\Omega\subset \mathbb R^n$ one has that
\begin{equation}\label{ani-Hardy}
{\displaystyle\int_{\Omega}  F_*(Du(x))^2{\rm d}x
	\geq \frac{(n-2)^2}{4}\int_{\Omega}\frac{u(x)^2}{F(x)^{2}}}{\rm d}x,\ \ \forall u\in W_0^{1,2}(\Omega);
\end{equation}
moreover, $\frac{(n-2)^2}{4}$ is optimal and never attained.

The following result is crucial in the study of extremal functions in the inequality $(\textbf{BPV})$.  
\begin{proposition}\label{lemma-swlsc-1} Let $(\mathbb R^n,F)$ be an $n$-dimensional reversible Minkowski space with the uniformity constant $l_F$ and fix  $\mu\in \left[0,l_F\frac{(n-2)^2}{4}\right]$. Then for every open set $\Omega\subset \mathbb R^n$, the functional 
	$$u\mapsto \mathcal K_\mu(u):= {\displaystyle\int_{\Omega}  F_*(Du(x))^2{\rm d}x
		-\mu\int_{\Omega}\frac{u(x)^2}{F(x)^{2}}}{\rm d}x $$ is positive and convex $($thus, sequentially weakly lower semicontinous$)$ on $W_0^{1,2}(\Omega).$
\end{proposition}

{\it Proof.} 
The positivity of $\mathcal{K}%
_\mu $ follows by the anisotropic Hardy inequality (\ref{ani-Hardy}) and $0<l_F\leq 1$. Let us fix $0<t<1$ and $u,v\in
W_{0}^{1,2}(\Omega)$. Then 
(\ref{eq:2uni})
and the anisotropic Hardy inequality (\ref{ani-Hardy}) imply
\begin{eqnarray*}
	\mathcal{K}_\mu\left(tu+(1-t)v\right) &=& \int_{\Omega}F_{*}(x,tD
	u(x)+(1-t)Dv(x))^2\mathrm{d}x -{\mu}\int_{\Omega}\frac{(tu(x)+(1-t)v(x))^2}{F(x)^2}\mathrm{d}x \\
	&\leq& t\int_{\Omega} F_*(x,Du(x))^2\mathrm{d}x
	+(1-t)\int_{\Omega}
	F_*(x,Dv(x))^2\mathrm{d}x \\
	&& -{l_{F}}t(1-t) \int_{\Omega} F_*(x,D(v-u)(x))^2\mathrm{d}x  -{\mu}\int_{\Omega}\frac{(tu(x)+(1-t)v(x))^2}{F(x)^2}\mathrm{d}x \\
	&=& t\mathcal{K}_\mu\left(u\right) +(1-t)\mathcal{K}_\mu\left(v\right)\\&& - t(1-t)l_F\int_{\Omega} \left( F_*(x,D(v-u)(x))^2-\mu l_F^{-1}\frac{%
		(v(x)-u(x))^2}{F(x)^2}\right)\mathrm{d}x \\
	&\leq& t\mathcal{K}_\mu\left(u\right) +(1-t)\mathcal{K}_\mu\left(v\right) ,
\end{eqnarray*}
which concludes the proof.
\hfill $\square$\\

\noindent {\bf {Proof of Theorem  \ref{theorem-extremals}}.} We split the proof into two parts. 

\underline{Case I}: $\alpha=0$ \textit{whenever} $n =2$, \textit{or} $\alpha>0$ \textit{whenever} $n \geq 3.$


\noindent Let $(\mathbb R^n,F)$ be a Minkowski space,  $\Omega\subset \mathbb R^n$ be an open bounded set  and $\alpha\in \left[\frac{n-2}{2}\sqrt{1-l_F^2},\frac{n-2}{2}\right]$ be fixed with the above properties; in particular, it turns out that $d_F(x_0,x)=F(x-x_0)$ for every $x_0,x\in \mathbb R^n$.   After translation, we may also consider that $x_0=0\in \Omega$.  If 
\begin{equation}\label{Reyligh}
\mu_{\alpha}(\Omega):=\inf_{u\in W_0^{1,2}(\Omega)}\left\{{\displaystyle\int_{\Omega}  F_*(Du(x))^2{\rm d}x
	-\left[\frac{(n-2)^2}{4}-\alpha^2\right]\int_{\Omega}\frac{u(x)^2}{F(x)^{2}}}{\rm d}x:\int_{\Omega}
{u(x)^2}{\rm d}x=1\right\},
\end{equation}
 it suffices to prove that 
\begin{equation}\label{amit-igazolni-kell}
\mu_{\alpha}(\Omega)\geq S_\alpha(\Omega).
\end{equation}
Moreover, since $F$ is absolutely homogeneous (so $F_*$), it is enough to consider only nonnegative test functions $u\in W_0^{1,2}(\Omega)_+$ in (\ref{Reyligh}).

Let us consider a minimizing sequence for $\mu_{\alpha}(\Omega)$, i.e., $\left\{u_k\right\}_k\subset W_0^{1,2}(\Omega)_+$ such that
\begin{equation}\label{limit-1}
{\int_{\Omega}
	{u_k(x)^2}}{\rm d}x=1, \ \ \forall k\in \mathbb N,
\end{equation}
and 
\begin{equation}\label{limit-2}
{\displaystyle\int_{\Omega}  F_*(Du_k(x))^2{\rm d}x
	-\left[\frac{(n-2)^2}{4}-\alpha^2\right]\int_{\Omega}\frac{u_k(x)^2}{F(x)^{2}}}{\rm d}x\to \mu_{\alpha}(\Omega)\ \ {\rm as}\ \ k\to \infty.
\end{equation} 
First, if $\alpha=0$ (when $n=2$),  by relation (\ref{limit-2}) we have that $\left\{u_k\right\}_k$ is bounded in $ W_0^{1,2}(\Omega)_+$. Second, if $\alpha>0$ (when $n\geq 3$), by the anisotropic Hardy inequality and (\ref{limit-2}), it follows again that $\left\{u_k\right\}_k$ is bounded in $ W_0^{1,2}(\Omega)_+$. Accordingly, there exists $\tilde u\in W_0^{1,2}(\Omega)_+$  such that (up to a subsequence) $\{u_k\}_k$ converges strongly to $\tilde u$ in $L^2(\Omega)$ and weakly to $\tilde u$ in $W_0^{1,2}(\Omega)$.   By (\ref{limit-1}) we directly have 
$${\int_{\Omega}
	{\tilde u(x)^2}}{\rm d}x=1,$$
and by Proposition \ref{lemma-swlsc-1} and relation (\ref{limit-2}),  

\begin{eqnarray*}
	\mu_{\alpha}(\Omega)&\leq& {\displaystyle\int_{\Omega}  F_*(D\tilde u(x))^2{\rm d}x
		-\left[\frac{(n-2)^2}{4}-\alpha^2\right]\int_{\Omega}\frac{\tilde u(x)^2}{F(x)^{2}}}{\rm d}x\\ &\leq& \liminf_{k\to \infty}\left[{\displaystyle\int_{\Omega}  F_*(Du_k(x))^2{\rm d}x
		-\left[\frac{(n-2)^2}{4}-\alpha^2\right]\int_{\Omega}\frac{u_k(x)^2}{F(x)^{2}}}{\rm d}x\right]\\
	&=&\mu_{\alpha}(\Omega). 
\end{eqnarray*}
Thus, $\tilde u\in W_0^{1,2}(\Omega)_+$ is a minimizer in (\ref{Reyligh}). Let $u^\star\in W_0^{1,2}(\Omega^\star)_+$ be the anisotropic symmetrization of $\tilde u$. 
By the anisotropic P\'olya-Szeg\H{o}, Hardy-Littlewood inequalities and Cavalieri principle one has
\begin{eqnarray}\label{PSZ-egyenloseg}
\mu_{\alpha}(\Omega)\nonumber&=& {\displaystyle\int_{\Omega}  F_*(D\tilde u(x))^2{\rm d}x
	-\left[\frac{(n-2)^2}{4}-\alpha^2\right]\int_{\Omega}\frac{\tilde u(x)^2}{F(x)^{2}}}{\rm d}x\\ &\geq& {\displaystyle\int_{\Omega^\star}  F_*(Du^\star(x))^2{\rm d}x
	-\left[\frac{(n-2)^2}{4}-\alpha^2\right]\int_{\Omega^\star}\frac{u^\star(x)^2}{F(x)^{2}}}{\rm d}x\\ \nonumber &\geq &
\inf_{v\in W_0^{1,2}(\Omega)_+}\left\{{\displaystyle\int_{\Omega^\star}  F_*(Dv^\star(x))^2{\rm d}x
	-\left[\frac{(n-2)^2}{4}-\alpha^2\right]\int_{\Omega^\star}\frac{v^\star(x)^2}{F(x)^{2}}}{\rm d}x:\int_{\Omega^\star}
{(v^\star)^2}{\rm d}x=1\right\}=:Q_{\alpha}(\Omega).
\end{eqnarray}

As above, one can  prove that the latter infimum is attained; let $v^\star\in W_0^{1,2}(\Omega^\star)_+$ be such a minimizer for $Q_{\alpha}(\Omega)$.   We may assume that
$v^\star\in C_0^1(\Omega^\star)_+$; otherwise, a density argument
applies. Thus, there exists a non-increasing function $h:[0,\infty)\to
[0,\infty)$ of class $C^1$ such that $v^\star(x)=h(\rho)$ where
$\rho=F(x).$ Since $\Omega^\star$ has a Wulff shape, there exists $R_\Omega>0$ such that $\omega_n R_\Omega^n=|\Omega^\star|=|\Omega|$; moreover, $h(R_\Omega)=0$. From the absolute homogeneity of $F$ and relation (\ref{tavolsag-derivalt}), we have that $$F_*(Dv^\star(x))=F_*(h'(\rho)D F(x))=-h'(\rho)F_*(D F(x))=-h'(\rho).$$ 
Consequently, the function $h$ is a minimizer for  
$$Q_{\alpha}(\Omega)=\min_{w\in C^1(0,R_\Omega);\ w(R_\Omega)=0}\frac{\displaystyle\int_0^{R_\Omega}(w'(\rho))^2r^{n-1}{\rm d}\rho-\left[\frac{(n-2)^2}{4}-\alpha^2\right]\displaystyle\int_0^{R_\Omega}w(\rho)^2\rho^{n-3}{\rm d}\rho}{\displaystyle\int_0^{R_\Omega}w(\rho)^2\rho^{n-1}{\rm d}\rho}.$$
The corresponding Euler-Lagrange equation for $h$ reads as 
\begin{equation}\label{diff-egy}
(h'(\rho)\rho^{n-1})'+\left[\frac{(n-2)^2}{4}-\alpha^2\right]h(\rho)\rho^{n-3}+Q_{\alpha}(\Omega)h(\rho)\rho^{n-1}=0, \ \ \rho\in (0,R_\Omega).
\end{equation}
If $w(\rho): = \rho^{\frac{n-2}{2}}h(\rho)$,  (\ref{diff-egy}) reduces to  the Bessel differential equation 
$$\rho^2w''(\rho)+\rho w'(\rho)+\left(Q_{\alpha}(\Omega)\rho^2-\alpha^2\right)w(\rho)=0.$$
Accordingly, (\ref{diff-egy}) has the general solution 
\begin{equation}\label{cand-extr}
h(\rho)=c_0\rho^\frac{2-n}{2}J_\alpha\left({\sqrt{Q_{\alpha}(\Omega)}}\rho\right)+c_1\rho^\frac{2-n}{2}Y_\alpha\left({\sqrt{Q_{\alpha}(\Omega)}}\rho\right), \ \ \rho\in (0,R_\Omega),
\end{equation}
where $c_0,c_1\in \mathbb R$, while $J_\alpha$ and $Y_\alpha$ are the Bessel functions of the first and second kind, respectively. Since $Y_\alpha$ is singular at the origin, we choose $c_1=0$; otherwise, $v^\star(x)=h(F(x))$ will not belong to $W_0^{1,2}(\Omega^\star)$. Furthermore, since $h(R_\Omega)=0$, it turns out that ${\sqrt{Q_{\alpha}(\Omega)}}R_\Omega=j_\alpha$, where $j_{\alpha}$ is the first
positive zero of $J_\alpha$, which gives that 
$${{Q_{\alpha}(\Omega)}}={{S_{\alpha}(\Omega)}} =j_\alpha^2\left(\frac{\omega_n}{|\Omega|}\right)^\frac{2}{n},$$ 
which implies (\ref{amit-igazolni-kell}), i.e, the validity of  {\rm $(\textbf{BPV})$} on $(\mathbb R^n,F)$.


If equality holds in  {\rm $(\textbf{BPV})$}, then we have that  $\mu_{\alpha}(\Omega)=S_{\alpha}(\Omega)$. The latter relation implies that in   (\ref{PSZ-egyenloseg}) we have equality; in particular, we have equality in the  the P\'olya-Szeg\H o inequality, i.e.,   
$$\displaystyle\int_{\Omega}  F_*(D\tilde u(x))^2{\rm d}x= \displaystyle\int_{\Omega^\star}  F_*(Du^\star(x))^2{\rm d}x.$$
Due to Esposito and Trombetti \cite[Theorem 5.1]{ET},  the latter relation implies that $\Omega=\Omega^*$ (up to translations) and $\tilde u$ agrees  almost everywhere (up to constant multiplication) with 
$$u^\star(x)=v^\star(x)=F(x)^\frac{2-n}{2}J_\alpha\left({\sqrt{S_{\alpha}(\Omega)}}F(x)\right),\ x\in \Omega^*,$$
see (\ref{cand-extr}). 
The asymptotic properties 
\begin{equation}\label{asymp}
J_\alpha(t)\sim \frac{1}{\Gamma(\alpha+1)}\left(\frac{t}{2}\right)^\alpha\ \ {\rm and}\ \ J'_\alpha(t)\sim \frac{1}{2\Gamma(\alpha)}\left(\frac{t}{2}\right)^{\alpha-1} \ {\rm  for}\ \  t\ll 1,
\end{equation}
 immediately imply that $u^\star\in W_0^{1,2}(\Omega^\star)$. 

Conversely, when $\Omega=\Omega^\star$, one clearly has $\mu_{\alpha}(\Omega)=S_{\alpha}(\Omega)$, and $u^\star$ is an extremal function in  {\rm $(\textbf{BPV})$}.   

\underline{Case II}: $\alpha=0$ \textit{whenever} $n\geq 3.$ In this case, we necessarily have that $l_F=1$, i.e., $F$ is Euclidean, thus we may proceed as in Brezis and V\'azquez \cite{BV}. Applying again (anisotropic) symmetrization, it is enough to prove (\textbf{BPV}) only for  symmetrized functions $v\in W_0^{1,2}(\Omega^\star).$ Let $v(x)=h(\rho)$ with $\rho=F(x),$ $x\in \Omega=\Omega^\star$. Note that   $h(R_\Omega)=0$,
where $R_\Omega=\left(\frac{|\Omega|}{\omega_n}\right)^\frac{1}{n}.$
Let
$w(\rho):=\rho^{\frac{n-2}{2}}h(\rho)$ with $\rho=F(x)$. Since $n\geq 3$, it turns out that $w(0)=0$. Moreover, since $w(R_\Omega)=0$, an integration by parts
gives
$$
\int_{\Omega^\star}  F_*(D v(x))^2{\rm d}x =n\omega_{n}\int_0^{R_\Omega}\rho (w'(\rho))^2{\rm d}\rho+n\omega_{n}\frac{(n-2)^2}{4}\int_0^{R_\Omega}\rho^{-1}w(\rho)^2{\rm d}\rho.
$$
Furthermore, one has
$$\int_{\Omega^\star}\frac{{v(x)}^2}{F(x)^{2}}{\rm d}x=n\omega_{n}\int_0^{R_\Omega}\rho^{-1}w(\rho)^2{\rm d}\rho$$  and $$ \int_{\Omega^\star}v(x)^2{\rm d}x=
n\omega_{n}\int_0^{R_\Omega}\rho w(\rho)^2{\rm d}\rho.$$
The above calculations show that in order to prove (\textbf{BPV}), it remains to check 
\begin{equation}\label{poin-2}
\int_0^{R_\Omega}\rho (w'(\rho))^2{\rm d}\rho\geq S_0(\Omega) \int_0^{R_\Omega}\rho w(\rho)^2{\rm d}\rho, 
\end{equation}
which is nothing but the optimal 2-dimensional Poincar\'e inequality. 

The equality in $(\textbf{BPV})$ would imply equality in (\ref{poin-2}). This would imply, similarly to Case I that $\Omega=\Omega^\star$ and $w(\rho)=J_0\left({\sqrt{S_{0}(\Omega)}}\rho\right).$ But  $v(x)=F(x)^{\frac{2-n}{2}}J_0\left({\sqrt{S_{0}(\Omega)}}F(x)\right)$ does not belong to $W_0^{1,2}(\Omega^\star)$, see (\ref{asymp}).  \hfill $\square$\\

\begin{remark}\rm A similar inequality to $(\textbf{BPV})$ can be stated also on not necessarily reversible Minkowski spaces. In such a setting, anisotropic symmetrization should be applied for positively homogeneous Minkowski norms, where  a set having a Wulff shape is homothetic to 
	the backward metric balls, see Van Schaftingen 	\cite{vanSch}. 
\end{remark}

\section{Proof of Theorems \ref{rigidity} and \ref{rigidity-Minkowski}}\label{section-4}

In order to provide the proof of Theorems \ref{rigidity} and \ref{rigidity-Minkowski}, we need some auxiliary results.

\begin{proposition}\label{layer-cake}
	Let $r>0$ and  $f:(0,r]\to \mathbb R$ be a non-increasing function such that $f(r)=0$ and $({M},F)$ be a complete $n$-dimensional reversible Finsler manifold $(n\geq 2)$ with nonnegative $n$-Ricci curvature. Then for every fixed $x_0\in M$, one has $$\int_{B_{x_0}(r)}f({d_F}(x_0,x)){\rm d}V_F(x)=\int_0^r {\sf A}_{x_0}(\rho)f(\rho){\rm d}\rho,$$
	where ${\sf A}_{x_0}(\rho):=\frac{\rm d}{{\rm d}\rho}{{\rm Vol}_F}(B_{x_0}(\rho))=\limsup_{\delta\to 0}\frac{{\rm Vol}_F(B_{x_0}(\rho+\delta))-{\rm Vol}_F(B_{x_0}(\rho))}{\delta}$ denotes the area of the sphere $\partial B_{x_0}(\rho)=\{y\in M:d_F(x_0,y)=\rho\}$.
\end{proposition}

{\it Proof.} By the  Bishop-Gromov volume comparison principle on $({M},F)$  we have that $\rho\mapsto \frac{{\rm Vol}_F(B_{x_0}(\rho))}{\rho^n}$ is non-increasing on $(0,\infty)$; in particular, $\rho\mapsto {{\rm Vol}_F(B_{x_0}(\rho))}$ is differentiable a.e.\ on $[0,\infty)$. Let $l_0=\lim_{\rho\to 0} f(\rho).$ By using the layer cake representation together with the facts that  $f:(0,r]\to \mathbb R$ is non-increasing and $f(r)=0$, an integration by parts gives
\begin{eqnarray*}
	\int_{B_{x_0}(r)}f({d_F}(x_0,x)){\rm d}V_F(x)&=&\int_0^{l_0}{\rm Vol}_F(\left\{x\in B_{x_0}(r):f({d_F}(x_0,x))>t\right\}){\rm d}t\\&=&
	\int_r^0{\rm Vol}_F( B_{x_0}(\rho))f'(\rho){\rm d}\rho \ \ \ \ \ \ \ \ [{\rm change\ of\ variables}\ t=f(\rho)]\\&=&-\int_0^r{\rm Vol}_F( B_{x_0}(\rho))f'(\rho){\rm d}\rho\\&=&\int_0^r \frac{\rm d}{{\rm d}\rho}{\rm Vol}_F(B_{x_0}(\rho))f(\rho){\rm d}\rho,
\end{eqnarray*}
which concludes the proof. 
\hfill $\square$\\



In the sequel we need some fine properties of  Bessel functions of the first kind; we first recall some basic properties of them which are well known in the literature, see e.g. Erd\'elyi,  Magnus,  Oberhettinger and Tricomi \cite{Bessel}. The Bessel function of the first kind with order $\alpha\in \mathbb R$ is the  solution of the differential equation
\begin{equation}\label{rec-23}
t^2y''(t)+ty'(t)+(t^2-\alpha^2)=0,
\end{equation}
which is nonsingular at the origin; we  denote it by $J_\alpha$. The function $J_\alpha$ has the following 
recurrence relations 
\begin{equation}\label{rec-0}
J_{\alpha+1}(t)+J_{\alpha-1}(t)=\frac{2\alpha}{t}J_\alpha(t), \ t>0;
\end{equation}
\begin{equation}\label{rec-1}
J'_\alpha(t)=-J_{\alpha+1}(t)+\frac{\alpha}{t}J_\alpha(t), \ t>0;
\end{equation}
\begin{equation}\label{rec-2}
J'_{\alpha+1}(t)=J_{\alpha}(t)-\frac{\alpha+1}{t}J_{\alpha+1}(t), \ t>0,
\end{equation}
where $J'_\alpha$ is the derivative of $J_\alpha$. We know that the positive zeros of $J_\alpha$ form an increasing sequence $\{j_{\alpha,k}\}_{k\in \mathbb N}$, and the Mittag-Leffler expansion yields  for every $\alpha>0$ that
\begin{equation}\label{Mittag-Lefler-1}
\frac{J_{\alpha+1}(t)}{J_{\alpha}(t)}=\sum_{k\geq1}\frac{2t}{j_{\alpha,k}^2-t^2},\ \ |t|<j_{\alpha,1};
\end{equation}
\begin{equation}\label{Mittag-Lefler-2}
t\frac{J'_{\alpha}(t)}{J_{\alpha}(t)}={\alpha}-\sum_{k\geq1}\frac{2t^2}{j_{\alpha,k}^2-t^2},\ \ |t|<j_{\alpha,1}.\end{equation}
In particular, by (\ref{Mittag-Lefler-1}) we easily obtain for $\alpha>-1$ the  Rayleigh sum
\begin{equation}\label{Bessel-zeros-sum}
\sum_{k\geq 1}\frac{1}{j_{\alpha,k}^2}=\frac{1}{4(\alpha+1)}.
\end{equation}
For simplicity, we use the  notation $j_\alpha:=j_{\alpha,1}$.  


\begin{proposition}\label{prop-monoton}
	Let $n\geq 2$  be an integer  and $\alpha\in \left[0,\frac{n-2}{2}\right]$. Then the following properties hold: 
	
	\begin{itemize}
		\item[(i)] for every $\beta \in [0,2]$,  the function $h_1({t}):={t}^{\beta-n}J_\alpha^2(j_\alpha {t})$ is non-increasing on $(0,1];$
		\item[(ii)] the function $h_2({t}):={t}^{1-n}J_\alpha(j_\alpha {t})J_{\alpha+1}(j_\alpha {t})$ is non-increasing on $(0,1];$
		\item[(iii)] the function $h_3({t}):={t}^{2-n}J_{\alpha+1}(j_\alpha {t})\left[J_{\alpha+1}(j_\alpha {t})-\frac{n+2\alpha}{j_\alpha {t}}J_\alpha(j_\alpha {t})\right]$ is non-decreasing on $(0,1].$
	\end{itemize}
\end{proposition}

{\it Proof.} (i) It is enough to prove that ${t}\mapsto \tilde h_1({t}):= \sqrt{h_1({t})}={t}^\frac{\beta-n}{2}J_\alpha(j_\alpha {t})$ is non-increasing on $(0,1).$ Since $J_\alpha(j_\alpha {t})\neq 0$ for ${t}\in (0,1)$, by relation (\ref{Mittag-Lefler-2}) 
we have 
$$\tilde h'_1({t})\frac{{t}^\frac{2+n-\beta}{2}}{J_\alpha(j_\alpha {t})}=\frac{\beta-n}{2}+j_\alpha {t}\frac{J'_{\alpha}(j_\alpha {t})}{J_{\alpha}(j_\alpha {t})}=\frac{\beta-n}{2}+\alpha -\sum_{k\geq1}\frac{2(j_\alpha {t})^2}{j_{\alpha,k}^2-(j_\alpha {t})^2}\leq \frac{\beta-n}{2}+\alpha\leq 0.$$

(ii) By using relations (\ref{rec-1}), (\ref{rec-2}), (\ref{Mittag-Lefler-1}) and (\ref{Bessel-zeros-sum}), we obtain for every ${t}\in (0,1)$ that 
\begin{eqnarray*}
	h'_2({t})\frac{{t}^n}{J_\alpha^2(j_\alpha {t})}&=&-j_\alpha {t}\frac{J^2_{\alpha+1}(j_\alpha {t})}{J^2_{\alpha}(j_\alpha {t})}-n\frac{J_{\alpha+1}(j_\alpha {t})}{J_{\alpha}(j_\alpha {t})}+j_\alpha {t}\\&\leq &-n\frac{J_{\alpha+1}(j_\alpha {t})}{J_{\alpha}(j_\alpha {t})}+j_\alpha {t}= j_\alpha {t}\left[-2n\sum_{k\geq1}\frac{1}{j_{\alpha,k}^2-(j_\alpha {t})^2}+1\right]\\ &\leq&  j_\alpha {t}\left[-2n\sum_{k\geq1}\frac{1}{j_{\alpha,k}^2}+1\right]=j_\alpha {t}\left[-\frac{n}{2(\alpha+1)}+1\right]\leq 0.
\end{eqnarray*}

(iii) A similar reasoning as in (ii) gives for every ${t}\in (0,1)$ that 
\begin{eqnarray*}
	h'_3({t})\frac{j_\alpha{t}^n}{J_\alpha^2(j_\alpha {t})[2(j_\alpha {t})^2+n(n+2\alpha)]}&= &\frac{J_{\alpha+1}(j_\alpha {t})}{J_{\alpha}(j_\alpha {t})}-\frac{(n+2\alpha)j_\alpha {t}}{2(j_\alpha {t})^2+n(n+2\alpha)} \\ &=&  j_\alpha {t}\left[2\sum_{k\geq1}\frac{1}{j_{\alpha,k}^2-(j_\alpha {t})^2}-\frac{(n+2\alpha)}{2(j_\alpha {t})^2+n(n+2\alpha)}\right]\\&\geq &j_\alpha {t}\left[\frac{1}{2(\alpha+1)}-\frac{1}{n}\right]\geq 0,
\end{eqnarray*}
which concludes the proof. 
\hfill $\square$\\

The following auxiliary result provides an unusual rigidity in the theory of functional inequalities  involving Bessel functions.  

\begin{proposition} \label{rigid-Bessel}	Let  $r>0$ be a real number, $n\geq 2$  be an integer, and $\alpha\in \left[0,\frac{n-2}{2}\right]$ be such  $\alpha>0$ whenever $n\geq 3.$
	Assume that a function $f:(0,r]\to [0,\infty)$ satisfies the following properties: 
	
	\begin{itemize}
		\item[(a)] $\liminf_{{t}\to 0}\frac{f({t})}{{t}^n}=1;$
		
		\item[(b)] the function ${t}\mapsto \frac{f({t})}{{t}^n}$ is non-increasing on $(0,r);$
		
		\item[(c)] $\displaystyle\int_0^1f(r{t}){t}^{-n+1}\left[J_{\alpha+1}^2(j_\alpha {t})-2\left(\frac{n-2}{2}-\alpha\right)\frac{J'_{\alpha}(j_\alpha {t})J_{\alpha}(j_\alpha {t})}{j_\alpha {t}}-J_{\alpha}^2(j_\alpha {t})\right]{\rm d}{t}\geq 0.$
	\end{itemize}
	Then $f({t})={t}^n$ for every ${t}\in (0,r).$
\end{proposition}

{\it Proof.} Two cases are distinguished, depending on $\alpha$ and $n.$  

\underline{Case I}: $\alpha>0$ \textit{whenever} $n\geq 3.$
For simplicity of notation, let us introduce the function $H_\alpha:(0,1]\to \mathbb R$ defined by
$$H_\alpha({t}):=J_{\alpha+1}^2(j_\alpha {t})-2\left(\frac{n-2}{2}-\alpha\right)\frac{J'_{\alpha}(j_\alpha {t})J_{\alpha}(j_\alpha {t})}{j_\alpha {t}}-J_{\alpha}^2(j_\alpha {t}).$$
By the integral identities   $$\int_0^1{t} J_{\alpha+1}^2(j_\alpha {t}){\rm d}{t}=\frac{J_{\alpha+1}^2(j_\alpha)}{2},\ \ \int_0^1 J_{\alpha}(j_\alpha {t})J_{\alpha}'(j_\alpha {t}){\rm d}{t}=0,\ \ \int_0^1{t} J_{\alpha}^2(j_\alpha {t}){\rm d}{t}=-\frac{J_{\alpha-1}(j_\alpha)J_{\alpha+1}(j_\alpha)}{2},$$
see formula (10.22.5) from \cite{Digital} and relation  (\ref{rec-0}),  we have that 
\begin{equation}\label{H-iden}
\displaystyle\int_0^1{t} H_\alpha({t}){\rm d}{t} =0.
\end{equation}

According to relation  (\ref{Mittag-Lefler-1}), the function ${t}\mapsto \frac{J^2_{\alpha+1}(j_\alpha {t})}{J^2_{\alpha}(j_\alpha {t})}$ is increasing on $(0,1)$ and one has 
$$\lim_{{t}\to 0}\frac{J^2_{\alpha+1}(j_\alpha {t})}{J^2_{\alpha}(j_\alpha {t})}=0\ \ {\rm and}\ \ \lim_{{t}\to 1}\frac{J^2_{\alpha+1}(j_\alpha {t})}{J^2_{\alpha}(j_\alpha {t})}=+\infty.$$
In a similar way, by (\ref{Mittag-Lefler-2}), the function ${t}\mapsto\frac{J'_{\alpha}(j_\alpha {t})}{j_\alpha {t} J_{\alpha}(j_\alpha {t})}$ is decreasing on $(0,1)$ and we have  $$\lim_{{t}\to 0}\frac{J'_{\alpha}(j_\alpha {t})}{j_\alpha {t} J_{\alpha}(j_\alpha {t})}=+\infty\ \ {\rm and}\ \ \lim_{{t}\to 1}\frac{J'_{\alpha}(j_\alpha {t})}{j_\alpha {t} J_{\alpha}(j_\alpha {t})}=-\infty.$$
The latter properties imply that the equation $H_\alpha({t})=0$ has a unique solution on $(0,1)$; let us denote by ${t}_\alpha^0\in (0,1)$ this element. The above analysis also shows that
\begin{equation}\label{rho-tulaj}
H_\alpha({t})<0,\ \forall {t}\in (0,{t}_\alpha^0)\ \ {\rm and}\ \  H_\alpha({t})>0,\ \forall {t}\in ({t}_\alpha^0,1].
\end{equation}

Due to (a) and (b), 
it follows that 
\begin{equation}\label{Funkc-rigid-1}
f({t})\leq {t}^n,\ \ \forall {t}\in (0,r).
\end{equation}
Relation (\ref{Funkc-rigid-1}) and the asymptotic properties (\ref{asymp}) of the Bessel function 
imply that the terms within inequality (c) are well defined. 
%
Let $g:(0,r]\to [0,1]$ be defined by $$g({t})=1-\frac{f({t})}{{t}^n}.$$
By (a) and (b), it follows that $\limsup_{{t}\to 0}{g({t})}=0$ and $g$ is non-decreasing on $(0,r).$ Since $f({t})={t}^n-g({t}){t}^n$, by means of (\ref{H-iden}), the inequality in (c) can be transformed equivalently into 
\begin{equation}\label{H-egyenlotlenseg}
\displaystyle\int_0^1 g(r{t}){t} H_\alpha({t}){\rm d}{t} \leq 0.
\end{equation}

We are going to prove that $g\equiv 0$ on $(0,r)$. 
To see this, we have 
\begin{eqnarray*}
	0&\geq& \displaystyle\int_0^1 g(r{t}){t} H_\alpha({t}){\rm d}{t}=\int_0^{{t}_\alpha^0} g(r{t}){t} H_\alpha({t}){\rm d}{t}+\int_{{t}_\alpha^0}^1 g(r{t}){t} H_\alpha({t}){\rm d}{t} \ \ \ \ \ \ [{\rm see}\ (\ref{H-egyenlotlenseg})] \\&\geq &g(r{t}_\alpha^0)\int_0^{{t}_\alpha^0} {t} H_\alpha({t}){\rm d}{t}+\int_{{t}_\alpha^0}^1 g(r{t}){t} H_\alpha({t}){\rm d}{t} \ \ \ \ \ \ [{\rm see}\ (\ref{rho-tulaj})\  {\rm and\ monotonicity\ of}\ g] \\&= &-g(r{t}_\alpha^0)\int_{{t}_\alpha^0}^1 {t} H_\alpha({t}){\rm d}{t}+\int_{{t}_\alpha^0}^1 g(r{t}){t} H_\alpha({t}){\rm d}{t} \ \ \ \ \ \ [{\rm see}\ (\ref{H-iden})] \\&= &\int_{{t}_\alpha^0}^1 \left[-g(r{t}_\alpha^0)+g(r{t})\right]{t} H_\alpha({t}){\rm d}{t}.
\end{eqnarray*}
By (\ref{rho-tulaj}) and the monotonicity of $g$ again, we necessarily have that $g(r{t})=g(r{t}_\alpha^0)$ for every ${t}\in ({t}_\alpha^0,1).$ Having this relation in  mind,  we have similarly as above that 
\begin{eqnarray*}
	0&\geq& \displaystyle\int_0^1 g(r{t}){t} H_\alpha({t}){\rm d}{t}=\int_0^{{t}_\alpha^0} g(r{t}){t} H_\alpha({t}){\rm d}{t}+g(r{t}_\alpha^0)\int_{{t}_\alpha^0}^1 {t} H_\alpha({t}){\rm d}{t}  \\&= &\int_0^{{t}_\alpha^0} g(r{t}){t} H_\alpha({t}){\rm d}{t}-g(r{t}_\alpha^0)\int_0^{{t}_\alpha^0} {t} H_\alpha({t}){\rm d}{t}  \\&= &\int_0^{{t}_\alpha^0} \left[g(r{t})-g(r{t}_\alpha^0)\right]{t} H_\alpha({t}){\rm d}{t}.
\end{eqnarray*}
Again, by (\ref{rho-tulaj}) and the monotonicity of $g$ we have $g(r{t})=g(r{t}_\alpha^0)$ for every ${t}\in (0,{t}_\alpha^0).$ Accordingly, $g(r{t})=g(r{t}_\alpha^0)$ for every ${t}\in (0,1).$ Since $\limsup_{{t}\to 0}{g({t})}=0$, we have that $g\equiv 0$ on $ (0,1),$ which concludes the proof in Case I.


\underline{Case II}: $\alpha=0$ \textit{whenever} $n =2.$ The proof is analogous to Case I; the only difference is that instead of $H_\alpha$ we consider the function 
$H_0:[0,1]\to \mathbb R$ defined by
$H_0({t}):=J_{1}^2(j_0 {t})-J_{0}^2(j_0 {t}).$
\hfill $\square$\\

We are now ready to prove  Theorem \ref{rigidity}. \\

\noindent \textbf{{Proof of Theorem \ref{rigidity}}.} We distinguish  two cases.   

\underline{Case I}: $\alpha>0$ \textit{whenever} $n\geq 3.$
Let us fix $r>0$ arbitrarily. We are going to prove that the function 
\begin{equation}\label{test-function}
u(x):=d_F(x_0,x)^\frac{2-n}{2}
J_\alpha
\left(j_\alpha
\frac{d_F(x_0,x)
}{r}
\right),\ x\in B_{x_0}(r),
\end{equation}
can be used as a test function in $(\textbf{BPV})$ on the open set $B_{x_0}(r)$. To do  this, we construct a sequence of functions with suitable  convergence properties. More precisely, for every $k\in \mathbb N$ with $k>\frac{1}{r}+1$, we consider the Lipschitz function
$u_{k}:B_{x_0}(r)\to \mathbb R$ defined by
\[ u_{k}(x)
=\left(\max\left\{\frac{1}{k},d_F(x_0,x)\right\}\right)^\frac{2-n}{2}
J_\alpha
\left(j_\alpha
\frac{\min\left\{\left(1+\frac{1}{k}\right)d_F(x_0,x),r\right\}
}{r}
\right). \]
Since $({M},F)$ is complete, the set
${\rm supp}(u_{k})=\overline {B_{x_0}(\frac{rk}{k+1})}$ is compact.
Therefore, by density reasons, $u_{k}$  can be uses as test functions
in $(\textbf{BPV})$ on $B_{x_0}(r)$, i.e., for every $k\in \mathbb N$ one has 
\begin{eqnarray}\label{limit-BPV}
\nonumber \int_{B_{x_0}(r)} F_*(x,Du_k(x))^2{\rm d}V_F(x) &\geq&
\left[\frac{(n-2)^2}{4}-\alpha^2\right]\int_{B_{x_0}(r)}\frac{u_k(x)^2}{d_F(x_0,x)^2}{\rm d}V_F(x)\\&&+S_{\alpha}(B_{x_0}(r))\int_{B_{x_0}(r)}
{u_k(x)^2}{\rm d}V_F(x).
\end{eqnarray}
Moreover, it immediately yields that 
\[  \lim_{k\to \infty}u_{k}(x)= u(x)\ \ {\rm and}\ \ \lim_{k\to \infty}F_*(x,Du_k(x))=F_*(x,Du(x))
\ \ {\rm a.e.}\ x\in B_{x_0}(r), \]
where 
\begin{equation}\label{weak-derivative}
F_*(x,Du(x))=d_F(x_0,x)^{-\frac{n}{2}}\left[\frac{n-2}{2} 
J_\alpha
\left(j_\alpha
\frac{d_F(x_0,x)
}{r}
\right)- j_\alpha
\frac{d_F(x_0,x)
}{r}J'_\alpha
\left(j_\alpha
\frac{d_F(x_0,x)
}{r}
\right)\right].
\end{equation}
In (\ref{weak-derivative}) we used the chain rule, relation (\ref{tavolsag-derivalt}) and the inequality $$j_\alpha{t}\frac{J'_\alpha(j_\alpha {t})}{J_\alpha(j_\alpha {t})}\leq \alpha \leq \frac{n-2}{2},\ \ {t}\in (0,1),$$
which follows by relation (\ref{Mittag-Lefler-2}).

We now prove that the functions 
$x\mapsto F_*(x,Du(x))^2$, $x\mapsto \frac{u(x)^2}{d_F(x_0,x)^2}$ and $x\mapsto u(x)^2$ belong to $L^1(B_{x_0}(r))$, respectively; these facts will be applied in the limiting process in (\ref{limit-BPV}) together with the Lebesgue dominated convergence theorem. 

First, by Proposition \ref{prop-monoton}/(i) (for $\beta=2$) and Proposition \ref{layer-cake} it turns out that
\begin{eqnarray}\label{kell-1}
\nonumber\int_{B_{x_0}(r)}
{u(x)^2}{\rm d}V_F(x)&=&\int_{B_{x_0}(r)}
d_F(x_0,x)^{2-n}
J^2_\alpha
\left(j_\alpha
\frac{d_F(x_0,x)
}{r}
\right){\rm d}V_F(x)\\&=&\nonumber r^{2-n}\int_{B_{x_0}(r)}
h_1\left(\frac{d_F(x_0,x)}{r}\right)
{\rm d}V_F(x)\\&=&\int_0^r {\sf A}_{x_0}(\rho)\rho^{2-n}J^2_\alpha
\left(j_\alpha
\frac{\rho
}{r}
\right){\rm d}\rho,
\end{eqnarray}
where $h_1$ comes from Proposition \ref{prop-monoton}. 
Since $\rho\mapsto \frac{{\rm Vol}_F(B_{x_0}(\rho))}{\rho^n}$ is  non-increasing, it follows (see also (\ref{volume-comp-altalanos-2})) that   
\begin{equation}\label{volume-compar}
{\rm Vol}_F(B_{x_0}(\rho))\leq \omega_n \rho^n\ \ {\rm and}\ \ {\sf A}_{x_0}(\rho)\leq n\omega_n\rho^{n-1},\ \ \forall \rho>0.
\end{equation}
Hence, by (\ref{asymp}) and (\ref{kell-1}) we have 
$$\int_0^r {\sf A}_{x_0}(\rho)\rho^{2-n}J^2_\alpha
\left(j_\alpha
\frac{\rho
}{r}
\right){\rm d}\rho\leq  n\omega_n \int_0^r \rho J^2_\alpha
\left(j_\alpha
\frac{\rho
}{r}
\right){\rm d}\rho
<+\infty.$$ 

Second, again by  Proposition \ref{prop-monoton}/(i) (for $\beta=0$) and Proposition \ref{layer-cake},  one has
\begin{eqnarray}\label{kell-2}
\nonumber\int_{B_{x_0}(r)}
\frac{u(x)^2}{d_F(x_0,x)^2}{\rm d}V_F(x)&=&\int_{B_{x_0}(r)}
d_F(x_0,x)^{-n}
J^2_\alpha
\left(j_\alpha
\frac{d_F(x_0,x)
}{r}
\right){\rm d}V_F(x)\\&=&\nonumber r^{-n}\int_{B_{x_0}(r)}
h_1\left(\frac{d_F(x_0,x)}{r}\right)
{\rm d}V_F(x)\\&=&\int_0^r {\sf A}_{x_0}(\rho)\rho^{-n}J^2_\alpha
\left(j_\alpha
\frac{\rho
}{r}
\right){\rm d}\rho\\&\leq &\nonumber
n\omega_n \int_0^r \rho^{-1} J^2_\alpha
\left(j_\alpha
\frac{\rho
}{r}
\right){\rm d}\rho<\infty,\nonumber
\end{eqnarray}
where we used (\ref{volume-compar}) and the asymptotic property 
(\ref{asymp}) together with the fact that $\alpha>0$. 

By using the recurrence relation (\ref{rec-1}) and the eikonal equation (\ref{tavolsag-derivalt}), it turns out that 
{\begin{equation}\label{ez-is-kelll}
F_*(x,Du(x))^2= d_F(x_0,x)^{-{n}}\left[\left(\frac{n-2}{2} -\alpha\right)
J_\alpha
\left(j_\alpha
\frac{d_F(x_0,x)
}{r}
\right)+ j_\alpha
\frac{d_F(x_0,x)
}{r}J_{\alpha+1}
\left(j_\alpha
\frac{d_F(x_0,x)
}{r}
\right)\right]^2.
\end{equation}}
Therefore, in order to handle the term $\displaystyle\int_{B_{x_0}(r)} F_*(x,Du(x))^2{\rm d}V_F(x)$, we have to study the behavior of three terms, coming from the development of (\ref{ez-is-kelll}). The first term appears precisely 
in (\ref{kell-2}). The second term is 
\begin{eqnarray}\label{kell-3}
\nonumber 0<I&:=&\int_{B_{x_0}(r)}
d_F(x_0,x)^{1-n}
J_\alpha
\left(j_\alpha
\frac{d_F(x_0,x)
}{r}
\right)J_{\alpha+1}
\left(j_\alpha
\frac{d_F(x_0,x)
}{r}
\right){\rm d}V_F(x)\\&=&\nonumber r^{1-n}\int_{B_{x_0}(r)}
h_2\left(\frac{d_F(x_0,x)}{r}\right)
{\rm d}V_F(x)\\&=&\int_0^r {\sf A}_{x_0}(\rho)\rho^{1-n}J_\alpha
\left(j_\alpha
\frac{\rho
}{r}
\right)J_{\alpha+1}
\left(j_\alpha
\frac{\rho
}{r}
\right){\rm d}\rho\\&\leq &\nonumber
n\omega_n \int_0^r  J_\alpha
\left(j_\alpha
\frac{\rho
}{r}
\right)J_{\alpha+1}
\left(j_\alpha
\frac{\rho
}{r}
\right){\rm d}\rho<\infty,\nonumber
\end{eqnarray}
where we used Proposition \ref{prop-monoton}/(ii), Proposition \ref{layer-cake} and relations (\ref{asymp}) and (\ref{volume-compar}). Finally, the third term is 
\begin{eqnarray}\label{kell-4}
\nonumber J&:=&\int_{B_{x_0}(r)}
d_F(x_0,x)^{2-n}
J_{\alpha+1}^2
\left(j_\alpha
\frac{d_F(x_0,x)
}{r}
\right){\rm d}V_F(x)\\&=&\nonumber r^{2-n}\int_{B_{x_0}(r)}
h_3\left(\frac{d_F(x_0,x)}{r}\right)
{\rm d}V_F(x)+\frac{r(n+2\alpha)}{j_\alpha}I.
\end{eqnarray}
Since $\rho\mapsto h_3(1)-h_3\left(\frac{\rho}{r}\right)$ is non-increasing on $(0,r]$, see Proposition \ref{prop-monoton}/(iii),  which vanishes at $r$, we may apply Proposition \ref{layer-cake} in order to obtain 
\begin{eqnarray*}
	\int_{B_{x_0}(r)}
	h_3\left(\frac{d_F(x_0,x)}{r}\right)
	{\rm d}V_F(x)&=&-\int_{B_{x_0}(r)}\left[h_3(1)-
	h_3\left(\frac{d_F(x_0,x)}{r}\right)\right]
	{\rm d}V_F(x)+h_3(1){\rm Vol}_F(B_{x_0}(r))\\&=&-\int_0^r {\sf A}_{x_0}(\rho)\left[h_3(1)-
	h_3\left(\frac{\rho}{r}\right)\right]{\rm d}\rho+h_3(1){\rm Vol}_F(B_{x_0}(r))\\&=&
	\int_0^r {\sf A}_{x_0}(\rho)
	h_3\left(\frac{\rho}{r}\right){\rm d}\rho.
\end{eqnarray*}
Accordingly, by the latter relation, (\ref{kell-3}) and (\ref{volume-compar}) one obtains 
\begin{eqnarray}\label{kell-5}
0<J&=&\nonumber r^{2-n}\int_0^r {\sf A}_{x_0}(\rho)
h_3\left(\frac{\rho}{r}\right){\rm d}\rho+\frac{r(n+2\alpha)}{j_\alpha}I\\&=&\int_0^r {\sf A}_{x_0}(\rho)\rho^{2-n}
J_{\alpha+1}^2\left(j_\alpha\frac{\rho}{r}\right){\rm d}\rho\\&\leq & \nonumber n\omega_n \int_0^r  \rho J_{\alpha+1}^2
\left(j_\alpha
\frac{\rho
}{r}
\right){\rm d}\rho<\infty.
\end{eqnarray}

Now, we are in the position to apply the Lebesgue dominated convergence theorem in (\ref{limit-BPV}), obtaining 
\begin{eqnarray}\label{limit-BPV-2}
\nonumber\int_{B_{x_0}(r)} F_*(x,Du(x))^2{\rm d}V_F(x) &\geq &
\left[\frac{(n-2)^2}{4}-\alpha^2\right]\int_{B_{x_0}(r)}\frac{u(x)^2}{d_F(x_0,x)^2}{\rm d}V_F(x)\\&&+S_{\alpha}(B_{x_0}(r))\int_{B_{x_0}(r)}
{u(x)^2}{\rm d}V_F(x).
\end{eqnarray}
By  (\ref{volume-compar})  we have 
\begin{equation}\label{masik-egyenlotlenseg}
S_{\alpha}(B_{x_0}(r))=j_{\alpha}^2\left(\frac{\omega_n}{{\rm Vol}_F(B_{x_0}(r))}\right)^\frac{2}{n}\geq \frac{j_{\alpha}^2}{r^2}.
\end{equation}
Therefore, by using the integral representations from (\ref{kell-1}),  (\ref{kell-2}),   (\ref{kell-3}) and (\ref{kell-5}), the inequalities (\ref{ez-is-kelll}), (\ref{limit-BPV-2}) and (\ref{masik-egyenlotlenseg}) imply after some elementary rearrangement that 
$$\int_0^r{\sf A}_{x_0}(\rho)\rho^{2-n}\left\{J^2_{\alpha+1}
\left(j_\alpha
\frac{\rho
}{r}
\right)-
2\left(\frac{n-2}{2}-\alpha\right)\frac{r}{j_\alpha \rho}\left[
-J_{\alpha+1}
\left(j_\alpha
\frac{\rho
}{r}\right)+\frac{\alpha r}{j_\alpha \rho}J_\alpha
\left(j_\alpha
\frac{\rho
}{r}
\right) \right]J_\alpha
\left(j_\alpha
\frac{\rho
}{r}
\right)-\right.$$$$\left.-J^2_{\alpha}
\left(j_\alpha
\frac{\rho
}{r}
\right)
\right\}{\rm d}\rho\geq 0.$$
We observe that for the second term in the above integrand we can apply the recurrence relation (\ref{rec-1}), obtaining  
$$\int_0^r{\sf A}_{x_0}(\rho)\rho^{2-n}\left[J^2_{\alpha+1}
\left(j_\alpha
\frac{\rho
}{r}
\right)-
2\left(\frac{n-2}{2}-\alpha\right)\frac{r}{j_\alpha \rho}J'_\alpha
\left(j_\alpha
\frac{\rho
}{r}
\right)J_\alpha
\left(j_\alpha
\frac{\rho
}{r}
\right)\right.\left.-J^2_{\alpha}
\left(j_\alpha
\frac{\rho
}{r}
\right)
\right]{\rm d}\rho\geq 0.$$
Now, we change the variable $\rho=r t,$ $t\in (0,1)$, which implies 
$$\int_0^1{\sf A}_{x_0}(rt)t^{2-n}\left[J^2_{\alpha+1}
\left(j_\alpha
t
\right)-
2\left(\frac{n-2}{2}-\alpha\right)\frac{1}{j_\alpha t}J'_\alpha
\left(j_\alpha
t
\right)J_\alpha
\left(j_\alpha
t
\right)-J^2_{\alpha}
\left(j_\alpha
t
\right)
\right]{\rm d}t\geq 0.$$
Since the latter inequality is valid for every $r>0$, an integration with respect to $r$ yields 
$$\int_0^1{\rm Vol}_F(B_{x_0}(rt))t^{1-n}\left[J^2_{\alpha+1}
\left(j_\alpha
t
\right)-
2\left(\frac{n-2}{2}-\alpha\right)\frac{1}{j_\alpha t}J'_\alpha
\left(j_\alpha
t
\right)J_\alpha
\left(j_\alpha
t
\right)-J^2_{\alpha}
\left(j_\alpha
t
\right)
\right]{\rm d}t\geq 0.$$

Due to (\ref{volume-comp-nullaban})  and the latter inequality, the assumptions in Proposition \ref{rigid-Bessel} are fulfilled by the non-increasing function $\rho\mapsto \frac{{\rm Vol}_F(B_{x_0}(\rho))}{\rho^n}$, $\rho\in (0,\infty)$. Therefore, it yields that $${\rm Vol}_F(B_{x_0}(\rho))=\omega_n\rho^n,\ \ \forall\rho>0.$$ 
We note that on $(M,F)$ with nonnegative $n$-Ricci curvature the latter relation does not depend on $x_0$, thus   
$${\rm Vol}_F(B_{x}(\rho))=\omega_n\rho^n,\ \ \forall x\in M,\rho>0,$$  
which implies that the flag curvature on $(M,F)$ is identically zero. \\

\underline{Case II}: $\alpha=0$ \textit{whenever} $n =2.$  The proof is similar to Case I. Let $r>0$ be  arbitrarily fixed. Instead of the function from (\ref{test-function}), we consider 
$$
u(x):=
J_0
\left(j_0
\frac{d_F(x_0,x)
}{r}
\right),\ x\in B_{x_0}(r).
$$
After a similar approximation procedure as above, we obtain by $(\textbf{BPV})$ that 
\begin{equation}\label{limit-BPV-3}
\int_{B_{x_0}(r)} F_*(x,Du(x))^2{\rm d}V_F(x) \geq
S_{0}(B_{x_0}(r))\int_{B_{x_0}(r)}
{u(x)^2}{\rm d}V_F(x).
\end{equation}
Since $$F_*(x,Du(x))=-\frac{j_0}{r}J_0'
\left(j_0
\frac{d_F(x_0,x)
}{r}\right)F_*\left(x,D d_F(x_0,\cdot)(x)\right)=\frac{j_0}{r}J_1
\left(j_0
\frac{d_F(x_0,x)
}{r}
\right)\ \ {\rm a.e.}\ x\in B_{x_0}(r),$$ 
and 
$$S_0(B_{x_0}(r))=j_{0}^2\frac{\omega_n}{{\rm Vol}_F(B_{x_0}(r))}\geq \frac{j_0^2}{r^2},$$
the inequality (\ref{limit-BPV-3}) implies 
$$\int_{B_{x_0}(r)} \left[J_1^2
\left(j_0
\frac{d_F(x_0,x)
}{r}
\right)-J_0^2
\left(j_0
\frac{d_F(x_0,x)
}{r}
\right)\right]{\rm d}V_F(x) \geq 0.$$
By the latter inequality we obtain 
$$\int_0^1{\rm Vol}_F(B_{x_0}(rt))t^{-1} \left[J_1^2
\left(j_0
t
\right)-J_0^2
\left(j_0
t
\right)\right]{\rm d}t \geq 0.$$	
It remains to apply Proposition \ref{rigid-Bessel} to obtain  ${\rm Vol}_F(B_{x_0}(\rho))=\omega_n\rho^n,$ $\rho>0.$ 
\hfill $\square$\\

\noindent {\bf {Proof of Theorem \ref{rigidity-Minkowski}}.} "(i)$\Rightarrow$(ii)" Since $(M,F)$ is Berwaldian, we have Ric$_n$=Ric, see Shen \cite{Shen-volume}, and we can apply Theorem \ref{rigidity}, obtaining that the flag curvature on $(M,F)$ is identically zero. Note that any Berwald space with vanishing flag curvature is necessarily a locally Minkowski space, see Bao, Chern and Shen \cite[Sect. 10.5]{BCS}. Now, the relation ${{\rm Vol}_F(B_x(\rho))}= \omega_n \rho^n$ for all $x\in M$ and $\rho>0$ implies that $(M,F)$ is in fact isometric to a Minkowski space. 

"(ii)$\Rightarrow$(i)" It follows directly by the first part of Theorem \ref{theorem-extremals}. \hfill $\square$\\

\section{Application in PDEs: Proof of Theorem \ref{alkalmazas}}\label{section-5}

\noindent \textbf{{Proof of Theorem \ref{alkalmazas}}}. Let us assume first that $({\mathcal P}_{\alpha,\lambda})$ has a nonzero solution, i.e., there exists $u\in W_0^{1,2}(B_0^F(1))_+\setminus \{0\}$ such that
\begin{equation}\label{diff-ujra-1}
-\Delta_F u(x)-\left[\frac{(n-2)^2}{4}-\alpha^2\right]\frac{u(x)}{F(x)^2}+\lambda u(x)  = u(x)^{p-1 },\ \  x\in B_0^F(1).
\end{equation}
Due to Theorem \ref{theorem-extremals}, the function 
$u^\star(x)=F(x)^\frac{2-n}{2}J_\alpha\left({{j_{\alpha}}}F(x)\right)$ belongs to  $W_0^{1,2}(B_0^F(1))$ and $u^\star(x)>0$ for every $x\in B_0^F(1)$. Moreover, the differential equation (\ref{diff-egy}) or a direct calculation yields that 
\begin{equation}\label{diff-ujra-2}
-\Delta_F u^\star(x)-\left[\frac{(n-2)^2}{4}-\alpha^2\right]\frac{u^\star(x)}{F(x)^2}=S_\alpha\left(B_0^F(1)\right)u^\star(x),\ \  x\in B_0^F(1).
\end{equation}
Note that  $S_\alpha\left(B_0^F(1)\right)=j_\alpha^2.$ By multiplying the equation (\ref{diff-ujra-1}) by $u^\star$, an integration by parts and (\ref{diff-ujra-2}) give that 
$$(\lambda+j_\alpha^2)\int_{B_0^F(1)}u^\star(x)u(x){\rm d}x=\int_{B_0^F(1)}u^\star(x)u(x)^{p-1}{\rm d}x>0,$$
which immediately yields $\lambda>-j_\alpha^2.$

Conversely, let us assume that $\lambda>-j_\alpha^2$ and define the positive number
$$c_{\alpha,\lambda}:=\left\{
\begin{array}{lll}
\min\left(1,1+\frac{\lambda}{j_0^2}\right),
& \hbox{if} &  n=2\ ({\rm and}\ \alpha=0), \\
\\
\frac{4\alpha^2}{(n-2)^2}\min\left(1,1+\frac{\lambda}{j_\alpha^2}\right), & \hbox{if} & n\geq 3\ ({\rm and}\ \alpha>0).
\end{array}\right.$$
Due to the validity of $(\textbf{BPV})$ (see Theorem \ref{theorem-extremals}), we have for every $u\in W_0^{1,2}(B_0^F(1))$ that
$$\mathcal K^2_{\alpha,\lambda}(u):=\int_{B_0^F(1)} \left\{ F_*(Du(x))^2
	-\left[\frac{(n-2)^2}{4}-\alpha^2\right]\frac{u(x)^2}{F(x)^{2}}+\lambda u(x)^2\right\}{\rm d}x\geq c_{\alpha,\lambda}\int_{B_0^F(1)}  F_*(Du(x))^2{\rm d}x.
$$
According to the classical Poincar\'e inequality and the latter inequality, $u\mapsto \mathcal K_{\alpha,\lambda}^{1/2}(u)$ defines a norm on $W_0^{1,2}(B_0^F(1))$, equivalent to the usual one.

Let $g,G:\mathbb R\to [0,\infty)$ be defined by $g(s)=s_+^{p-1}$ and $G(s)=\frac{s_+^{p}}{p},$ where $s_+=\max(s,0).$ We associate with problem $({\mathcal P}_{\alpha,\lambda})$ its energy functional $\mathcal E_{\alpha,\lambda}:W_0^{1,2}(B_0^F(1))\to \mathbb R$ defined by
$$\mathcal E_{\alpha,\lambda}(u)=\frac{1}{2}\mathcal K^2_{\alpha,\lambda}(u)-\int_{B_0^F(1)}G(u(x)){\rm d}x.$$
In a standard manner we can prove that $\mathcal E_{\alpha,\lambda}\in C^1(W_0^{1,2}(B_0^F(1)); \mathbb R)$. Moreover, since $W_0^{1,2}(B_0^F(1))$ can be compactly embedded into $L^p(B_0^F(1))$ ($p\in (2,2^*)$) and $s\mapsto g(s)$ verifies the usual Ambrosetti-Rabinowitz condition (see e.g. Willem \cite[Lemma 1.20]{Willem}), it turns out that $\mathcal E_{\alpha,\lambda}$ satisfies the Palais-Smale condition at each level. Moreover, since $p>2$, $\mathcal E_{\alpha,\lambda}$ satisfies also the mountain pass geometry; namely, 
\begin{itemize}
	\item there exists a sufficiently small  $\rho>0$ such that
	$$\inf_{\mathcal K_{\alpha,\lambda}(u)=\rho}\mathcal E_{\alpha,\lambda}(u)>0=\mathcal E_{\alpha,\lambda}(0);$$
	\item for sufficiently large $t>0$  and $u^\star$ coming from (\ref{diff-ujra-2}) one has 
	$$\mathcal E_{\alpha,\lambda}(tu^\star)=t^2\mathcal K^2_{\alpha,\lambda}(u^\star)-t^p\int_{B_0^F(1)}G(u^\star(x)){\rm d}x<0.$$

\end{itemize}
Therefore, by the mountain pass theorem it follows the existence of a critical point $u\in W_0^{1,2}(B_0^F(1))$ for $\mathcal E_{\alpha,\lambda}$ with positive energy level (thus $u\neq 0$), which is precisely a weak solution to problem  
\[ \   \left\{ \begin{array}{lll}
-\Delta_F u(x)-\left[\frac{(n-2)^2}{4}-\alpha^2\right]\frac{u(x)}{F(x)^2}+\lambda u(x)  = (u(x))_+^{p-1 },& &  x\in B_0^F(1); \\
 u\in W_0^{1,2}(B_0^F(1)).
\end{array}\right. \]
Multiplying the latter equation by $u_-(x)=\min(u(x),0)$, an integration on $B_0^F(1)$ gives that 
$\mathcal K^2_{\alpha,\lambda}(u_-)=0$, which implies $u_-=0$. Therefore, $u\geq 0$ is a solution for the original problem $({\mathcal P}_{\alpha,\lambda})$, which concludes the proof. \hfill $\square$\\

\noindent {\bf Acknowledgment.} The authors thank  \'Arp\'ad Baricz for useful conversations on Bessel functions.


\begin{thebibliography}{99}

\bibitem{ACR} Adimurthi, N. Chaudhuri, M. Ramaswamy, An improved Hardy-Sobolev
inequality and its applications, Proc. Amer. Math. Soc. 130 (2002),
489--505.

\bibitem{AIHP-Lions} A. Alvino, V. Ferone, P.-L. Lions, G.
Trombetti, Convex symmetrization and applications,  Ann. Inst. H.
Poincar\'e Anal. Non Lin\'eaire  14 (1997),  no. 2, 275--293.



\bibitem{Barbatis} G. Barbatis, S. Filippas, A. Tertikas,
A unified approach to improved $L^p$ Hardy inequalities with best
constants, Trans. Amer. Math. Soc. 356 (2004), 2169--2196.

\bibitem{BCS} D.~Bao, S.~S.~Chern, Z.~Shen, {Introduction to
Riemann--Finsler Geometry,} Graduate Texts in Mathematics, 200,
Springer Verlag, 2000.

%

\bibitem{Grillo} E. Berchio, D. Ganguly, G. Grillo,  Sharp Poincar\'e-Hardy and Poincar\'e-Rellich inequalities on the hyperbolic space, J. Funct. Anal. 272 (2017), no. 4, 1661--1703.


\bibitem{BV} H. Brezis, J. L. V\'azquez, Blowup solutions of some nonlinear elliptic
problems, Revista Mat. Univ. Complutense Madrid 10 (1997), 443--469.

\bibitem{Caron} {G.\  Carron,} In\'egalit\'es de Hardy sur les vari\'et\'es riemanniennes non-compactes, 
{J. Math. Pures Appl.} (9) 76 (1997), no. 10, 883--891.

\bibitem{Carron-2} {G.\  Carron,} In\'egalit\'es isop\'erim\'etriques de Faber-Krahn et consequences, Publications de
l'Institut Fourier, 220, 1992.


\bibitem{Cheng} S.-Y. Cheng, Eigenvalue comparison theorems
and its geometric applications, Math. Z. 143 (1975), 289--297. 


\bibitem{D-D-olaszok} {L.\  D'Ambrosio,
	S.\  Dipierro,} Hardy inequalities on Riemannian manifolds and
applications, Ann. Inst. H.
Poincar\'e Anal. Non Lin\'eaire,  31 (2014), no. 3, 449--475.





\bibitem{Bessel} A. Erd\'elyi, W. Magnus, F. Oberhettinger, F. Tricomi,  Higher transcendental functions. Vols. I, II. Based, in part, on notes left by Harry Bateman. McGraw-Hill Book Company, Inc., New York-Toronto-London, 1953. 

\bibitem{ET} L. Esposito, C. Trombetti, 
Convex symmetrization and P\'olya-Szeg\H o inequality,
Nonlinear Anal. 56 (2004), no. 1, 43--62. 

\bibitem{FKV} {C.\ Farkas, A.\  Krist\'aly, C.\ Varga}, Singular Poisson equations on Finsler-Hadamard manifolds, {Calc. Var. Partial Differential Equations}  54 (2015), no. 2, 1219--1241. 


\bibitem{GM} N. Ghoussoub,  A. Moradifam, Bessel pairs and optimal Hardy
and Hardy-Rellich inequalities, Math. Ann. 349 (2011), no. 1, 1--57.

\bibitem{GM-2} N. Ghoussoub,  A. Moradifam, On the best possible remaining term in the Hardy inequality, Proc. Natl. Acad. Sci. USA 105 (2008), no. 37, 13746--13751.


\bibitem{KO-1} {I.\  Kombe, M.\  \"Ozaydin,} Improved Hardy and Rellich
inequalities on Riemannian manifolds, Trans. Amer. Math. Soc. 361
(2009), no. 12, 6191--6203.

\bibitem{KO-2} {I.\  Kombe, M.\  \"Ozaydin,} Hardy-Poincar\'e, Rellich and uncertainty principle inequalities on Riemannian manifolds,  Trans. Amer. Math. Soc. 365 (2013), no. 10,
5035--5050.


\bibitem{Kristaly-JMPA} {A.\  Krist\'aly}, Sharp uncertainty principles on Riemannian manifolds: the influence of curvature, J. Math. Pures Appl., in press. DOI: 10.1016/j.matpur.2017.09.002.

\bibitem{Kri-Ohta} {A.\  Krist\'aly, S.\  Ohta,} Caffarelli-Kohn-Nirenberg inequality on metric measure spaces with applications, {Math. Ann.} 357 (2013), no. 2, 711--726.



\bibitem{Digital} F. W. J. Olver, D. W. Lozier, R. F. Boisvert, C. W. Clark (eds.), NIST Handbook
of Mathematical Functions, Cambridge University Press, Cambridge, 2010.



\bibitem{Ohta-Calculus} S. Ohta, Finsler interpolation inequalities, Calc. Var. Partial Differential Equations 36 (2009), no. 2, 211--249.


 
 \bibitem{Ohta-Sturm}  {S.\  Ohta, K.-T.\  Sturm,} Heat flow on Finsler manifolds. {Comm. Pure and Appl. Math.}  62 (2009), no. 10,
 1386--1433.
 
 \bibitem{Shen-volume} Z. Shen, Volume comparison and its applications in Riemann-Finsler geometry,
 Adv. Math. {128} (1997), no. 2, 306--328.
 
%
%
 

%


\bibitem{vanSch} J. Van Schaftingen, Anisotropic symmetrization,  Ann. Inst. H.
Poincar\'e Anal. Non Lin\'eaire, 23 (2006),  no. 4, 539-565.



\bibitem{Willem} M. Willem, Minimax theorems.
Progress in Nonlinear Differential Equations and their Applications, 24. Birkh\"auser Boston, Inc., Boston, MA, 1996.




\bibitem{YSK} Q. Yang, D. Su, Y. Kong, Hardy inequalities on Riemannian manifolds with negative curvature, Commun. Contemp. Math. 16(2) (2014), Article ID: 1350043, 24 pp.

\end{thebibliography}
\end{document}